\documentclass[a4paper,12pt]{article}
\usepackage{latexsym}

\usepackage{latexcad}

\newcommand{\pic}[2]{\BoxedEPSF{#1 scaled #2}}

\usepackage[hideboxes]{boxedeps}

\newcommand\Doublehex{\pic{Doublehex.ART}{400}}
\newcommand\Doublehexknot{\pic{Doublehexknot.ART}{400}}
\newcommand\Doublehexcircle{\pic{Doublehexcircle.ART}{400}}

\newcommand\Splittings{\pic{Splittings.ART}{400}}

\newcommand\Chord{\pic{Chord.ART}{400}}
\newcommand\Ychord{\pic{Ychord.ART}{400}}
\newcommand\Ytwochord{\pic{Ytwochord.ART}{400}}

\newcommand\figtwo{\pic{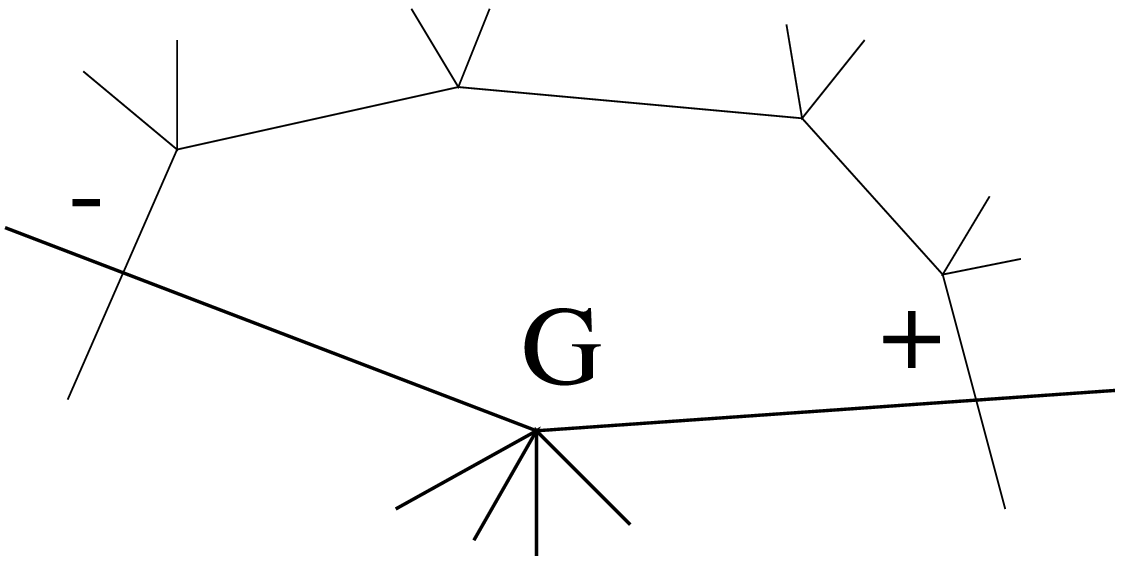}{400}}
\newcommand\figthreea{\pic{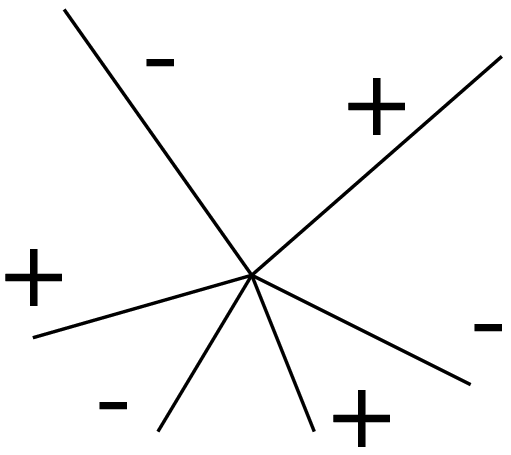}{300}}
\newcommand\figthreeb{\pic{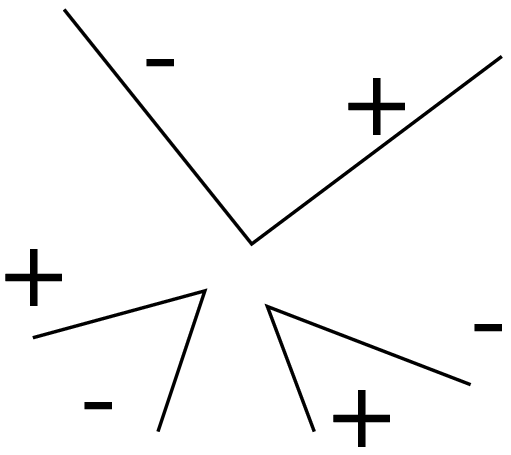}{300}}
\newcommand\figfour{\pic{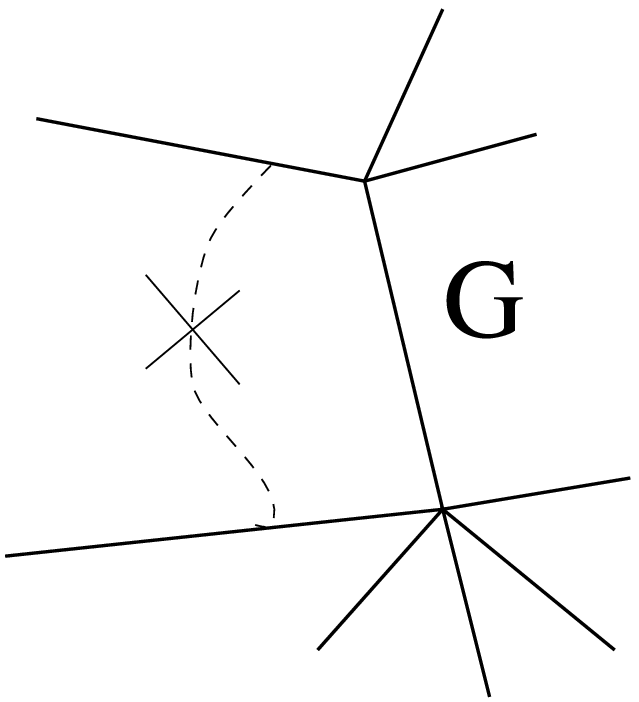}{400}}
\newcommand\figsevena{\pic{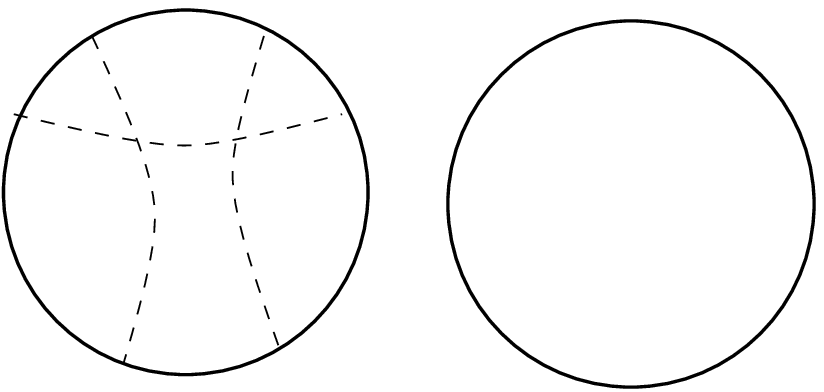}{300}}
\newcommand\figsevenb{\pic{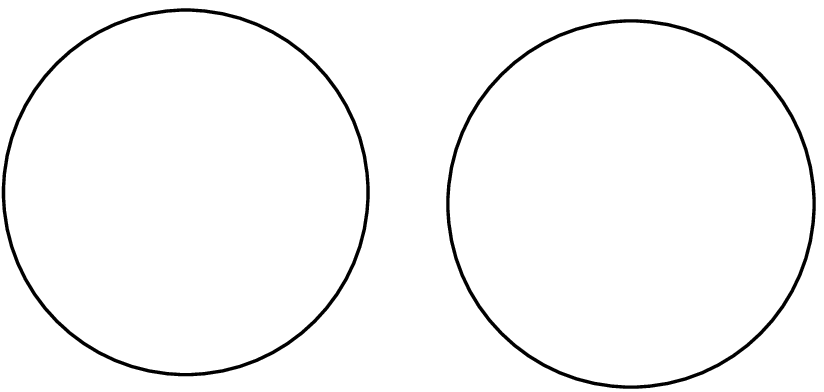}{300}}

\newcommand\figfive{\pic{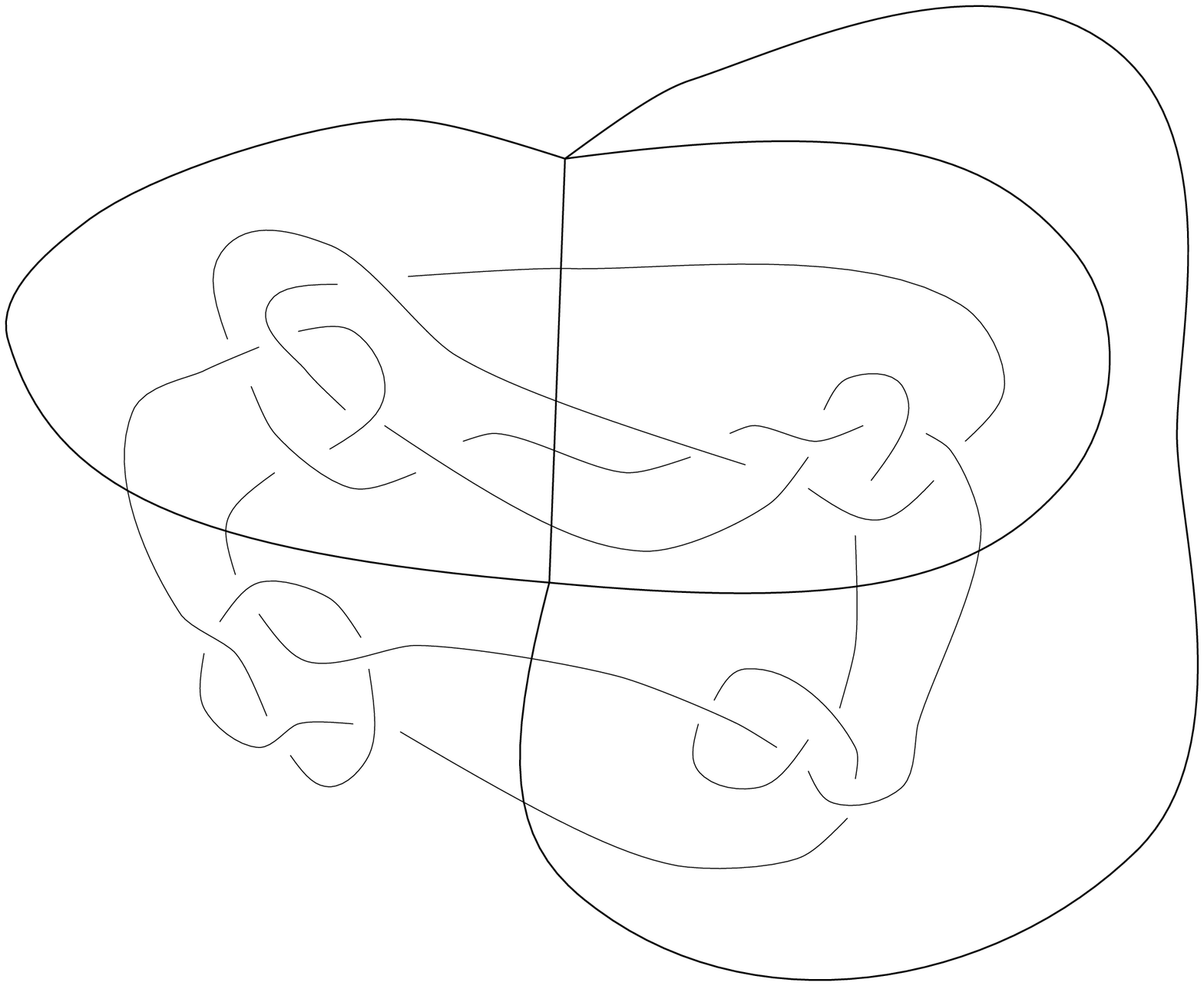}{300}}
\newcommand\figfivega{\pic{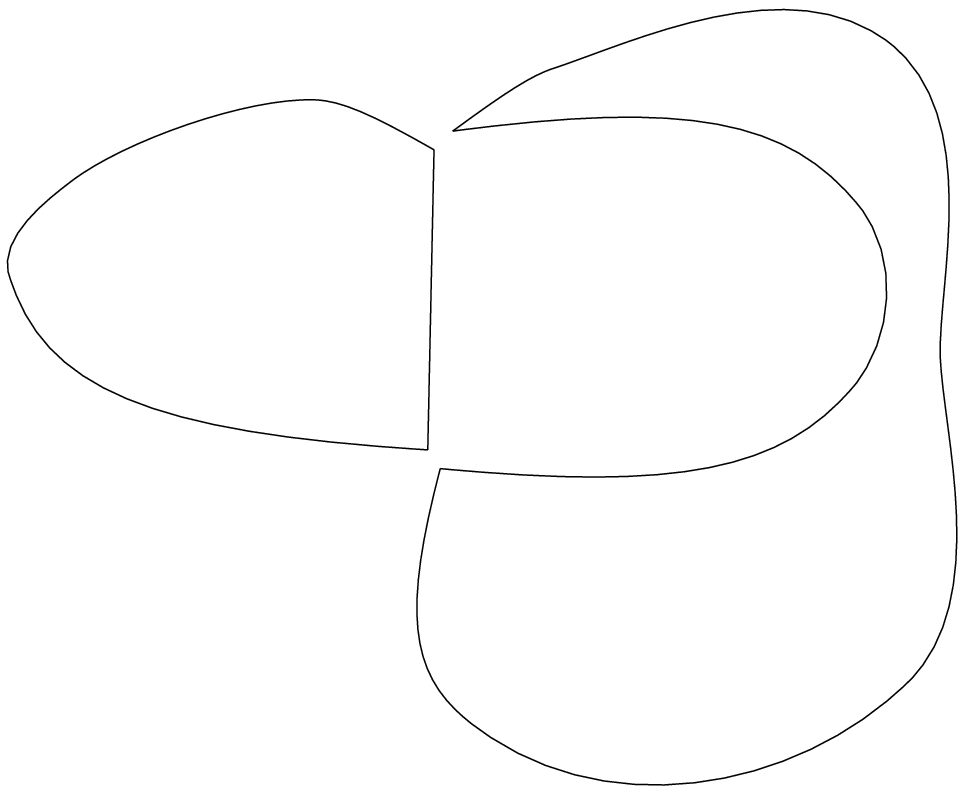}{300}}
\newcommand\figfivegb{\pic{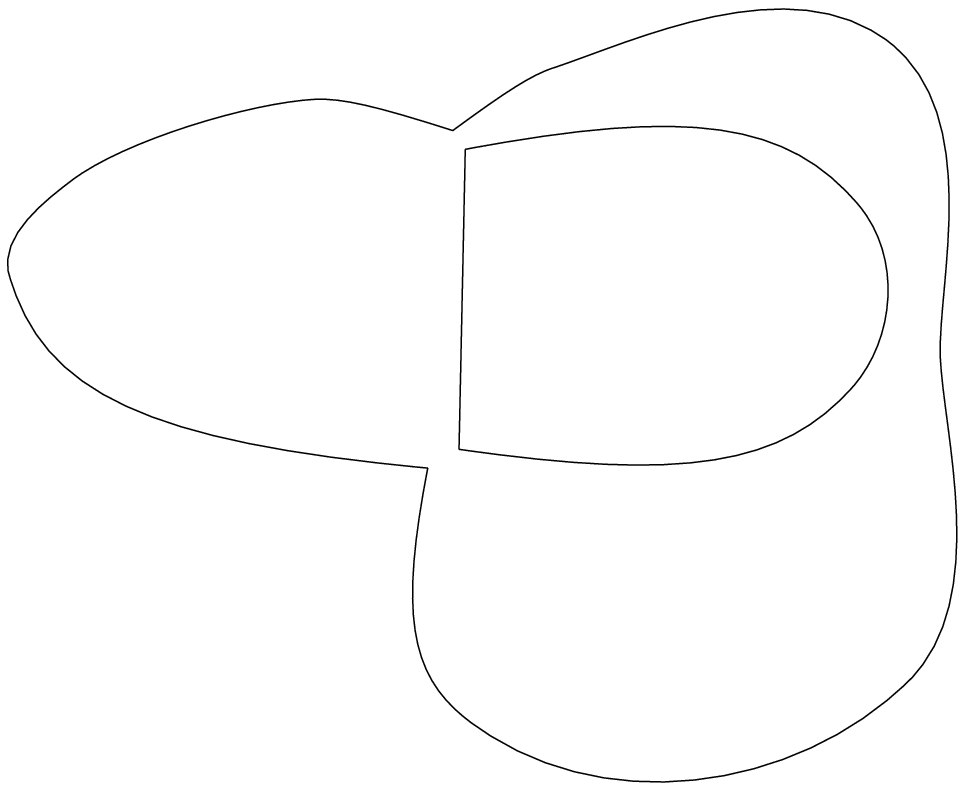}{300}}
\newcommand\figsixa{\pic{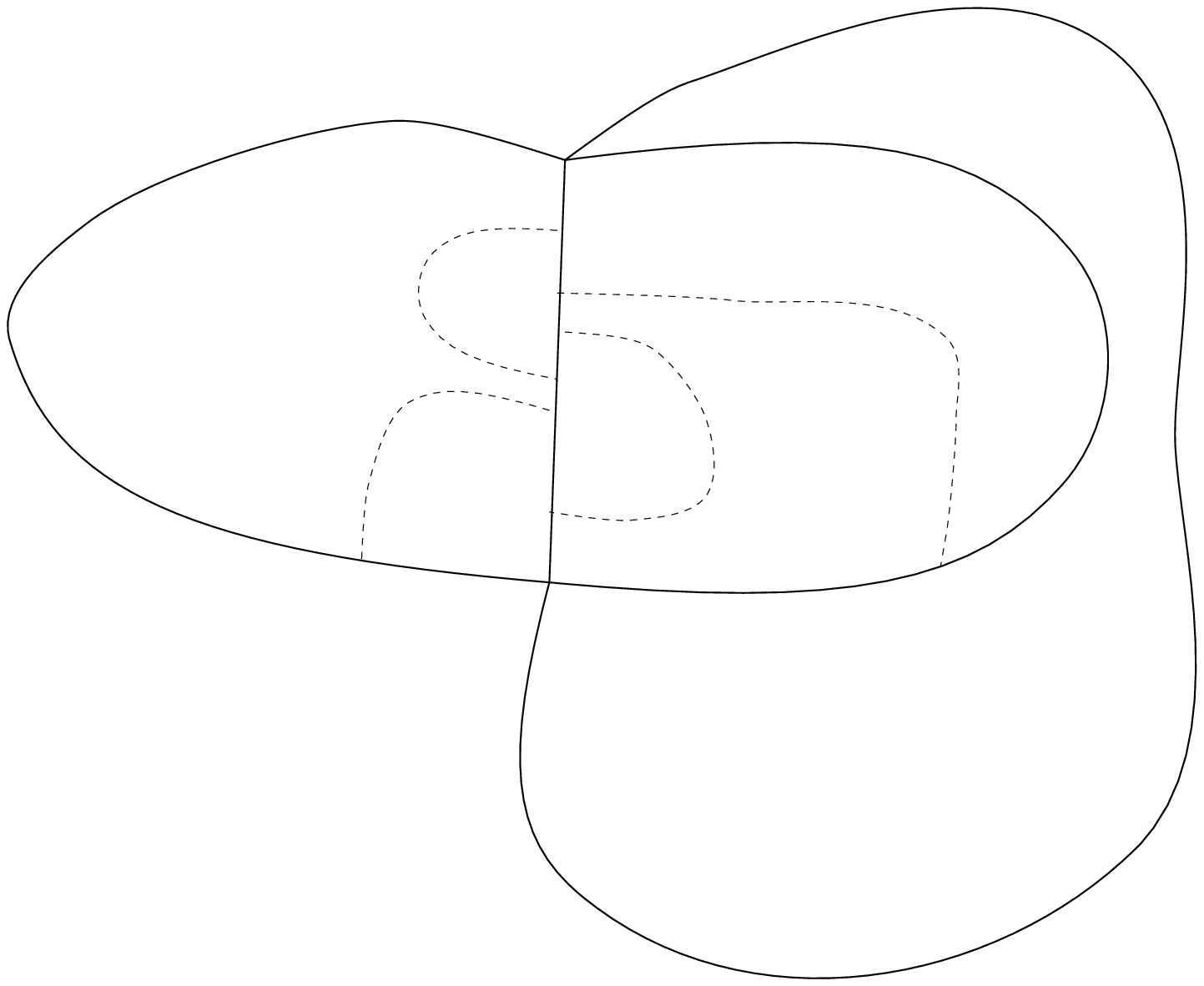}{250}}
\newcommand\figsixb{\pic{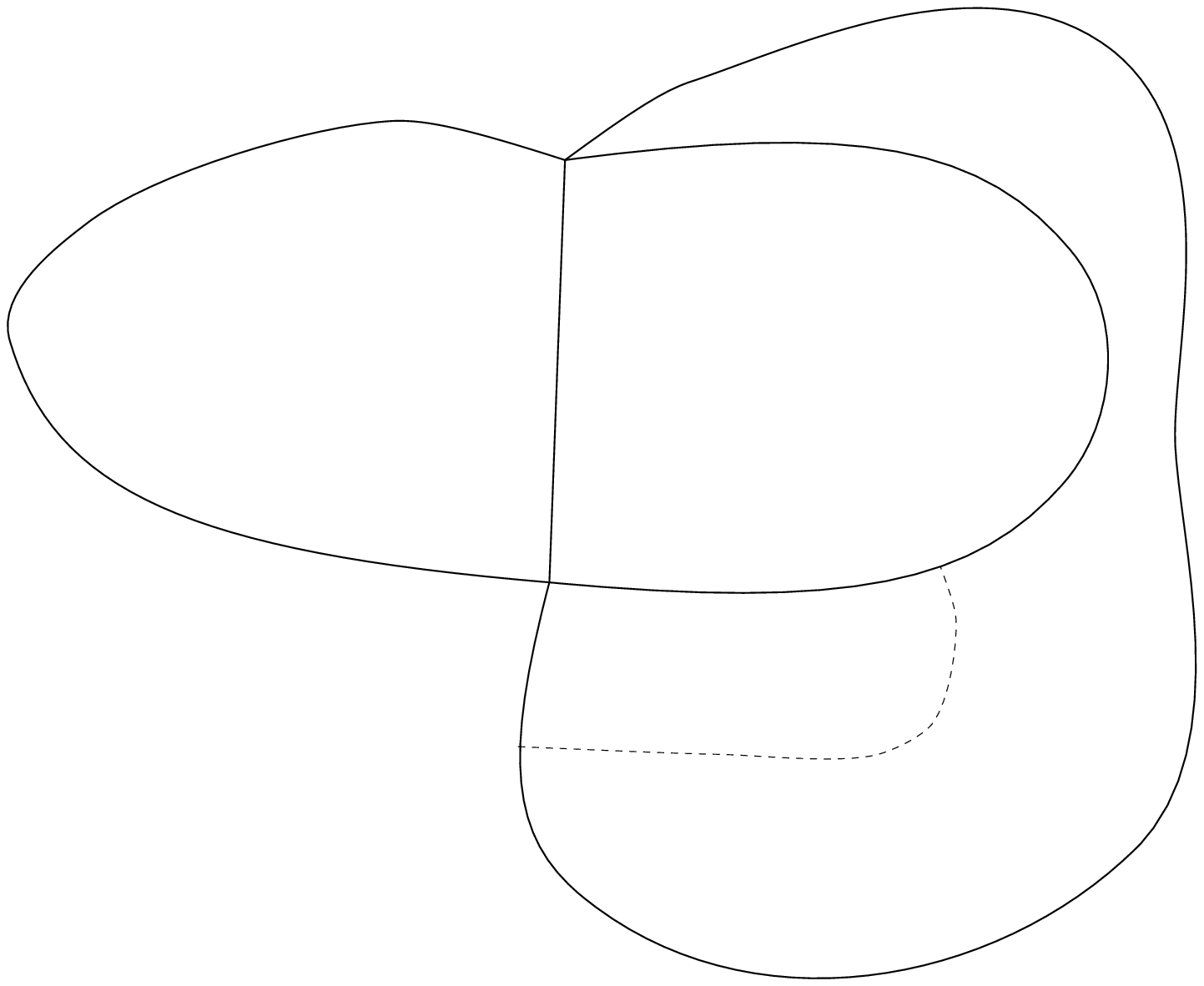}{250}}
\newcommand\figeight{\pic{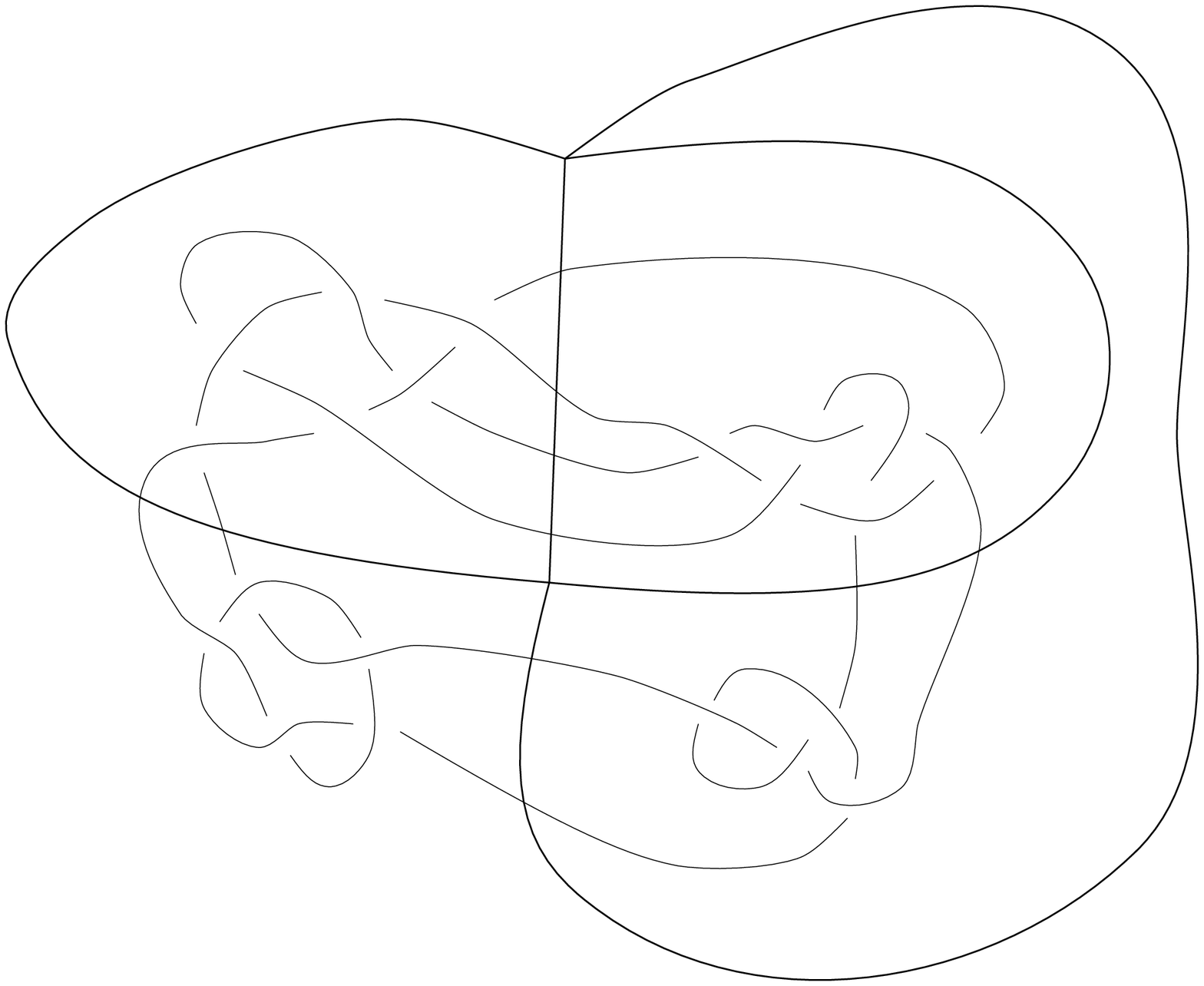}{300}}
\newcommand\fignine{\pic{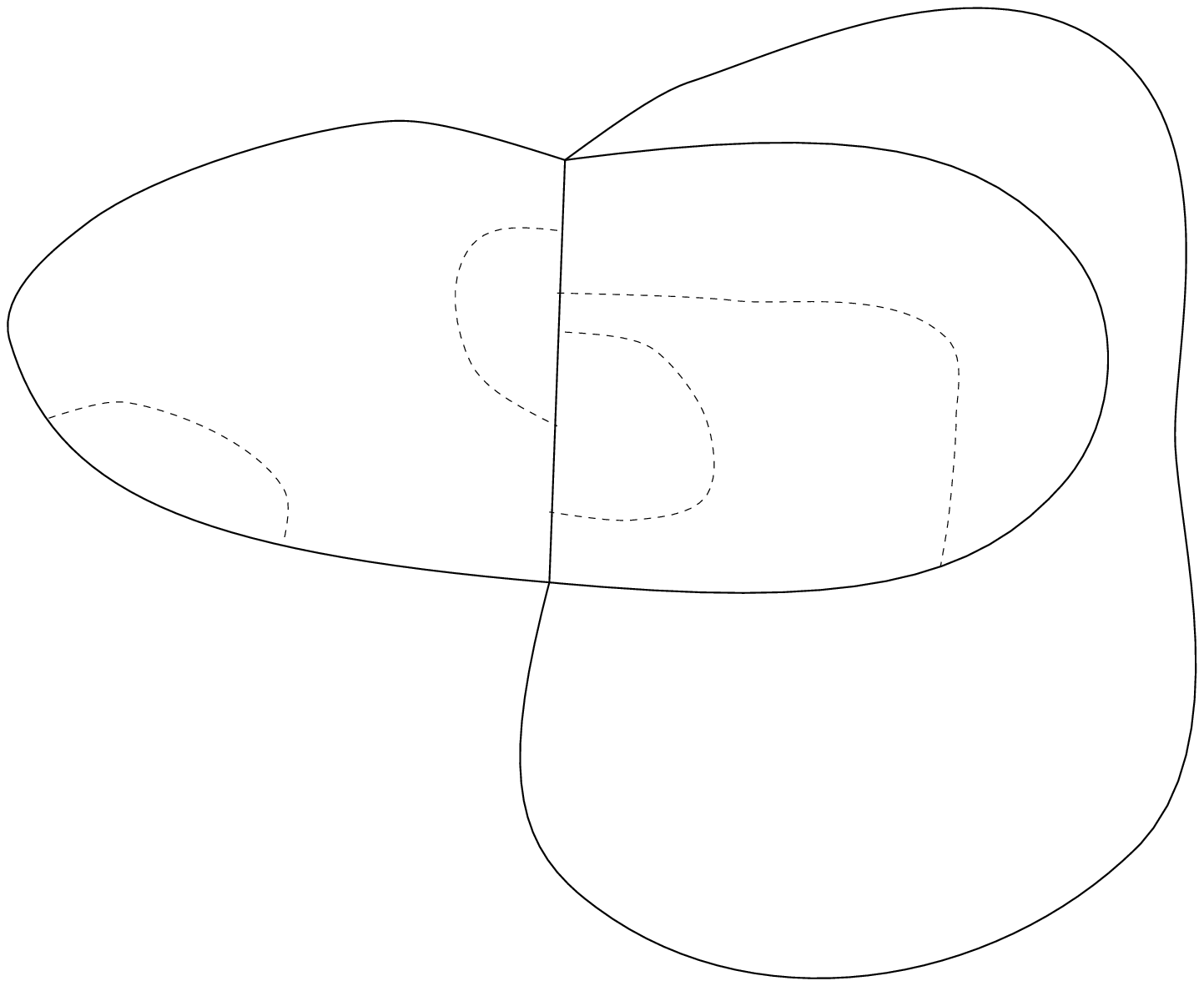}{300}}
\newcommand\figten{\pic{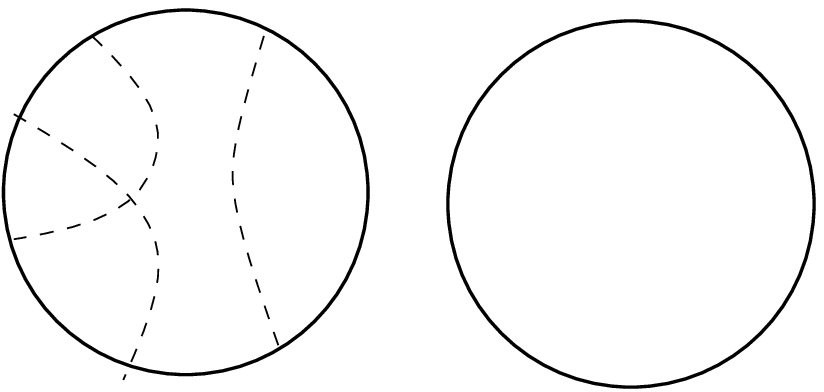}{300}}

\newcommand{\abs}[1]{\left\vert#1\right\vert}
\newcommand{\seq}[1]{\left<#1\right>}
\newcommand{\brac}[1]{\langle #1\rangle}
\newcommand{\bd}{\brac D}
\newtheorem{theorem}{Theorem}

\newtheorem{property}{Property}

\newenvironment{remark}{\par\smallskip%
\noindent\textbf{Remark.}\  }%
{\par\smallskip}

\newenvironment{Def}{\par\smallskip%
\noindent\textbf{Definition.}\  }%
{\par\smallskip}

\newenvironment{Proof}[1][{}]{\par\smallskip%
\noindent\textit{Proof #1: }\  }
{\hfill$\Box$\par\smallskip}

\begin{document}
\title{The spread and extreme terms of Jones polynomials}
\author{Yongju Bae and H.R.Morton}

\maketitle
\begin{abstract} We adapt Thistlethwaite's alternating tangle
decomposition of a knot diagram to identify the potential extreme
terms in its bracket polynomial, and give a simple combinatorial
calculation for their coefficients, based on the intersection graph
of certain chord diagrams.
\end{abstract}

\section*{Introduction}

One of the most striking combinatorial applications of the Jones
polynomial has been the result of Murasugi and Thistlethwaite which
characterises alternating knots by relating the spread of the Jones
polynomial to the number of crossings in the knot diagram. The result
follows from the identification of the two potential extreme terms in
the bracket polynomial and the calculation that each term occurs with
coefficient $\pm1$.

Lickorish and Thistlethwaite \cite{Lick} were able to widen the class
of knots for which a similar exact bound for the spread of the
polynomial could be given. In this wider class of `adequate' diagrams
they were again able to specify extreme terms and show that their
coefficients were
$\pm1$. Later Thistlethwaite \cite{Th} made use of a natural
decomposition of a given link diagram into maximal alternating
pieces, and formulated a bound for the spread purely in terms of
combinatorial data from the non-alternating part of the diagram. 

In this paper we extend his approach, again looking at the structure
of a given diagram in terms of its maximal alternating tangles. We
incorporate some `boundary connection' information about each
alternating tangle, which depends only on the immediate neighbourhood
of the boundary of the tangle, and combine this with the data from
the non-alternating part to give a simple combinatorial calculation
for the coefficients of two potential extreme terms in the bracket
polynomial for the diagram. 

In the case of adequate diagrams our calculations
immediately give $\pm1$ for these coefficients, and hence the result
on the spread of such diagrams.  In more general conditions the
extreme terms considered will depend on the diagram chosen. When our
extreme coefficients are  both non-zero, we can again identify the
spread of the Jones polynomial, while if one of the coefficients is
zero we then at least get a better upper bound  for the
spread than the initial one from the given diagram.

\section{The extreme states bound} We recall the states sum description
of the Kauffman bracket polynomial $\bd$ for an unoriented link
diagram $D$.

Label the four quadrants at each crossing $A$ or
$B$, according to the rule that the overcrossing strand sweeps out
the $A$ quadrants when turned anticlockwise. The two possible local
splittings  are termed  the
$A$-split and the $B$-split.

A {\em state} $s$ of $D$ is a labelling of each crossing  by either
$A$ or $B$. Making the corresponding split for each crossing gives a
number of disjoint embedded closed curves, called the {\em state
circles} for $s$. We write
$|s|$ for the number of state circles, and $a(s)$, $b(s)$ respectively
for the number of crossings labelled
$A$, $B$ in the state $s$. Following Kauffman \cite{K} we retain
information about the original crossings, in the form of chords
on the split diagram, which we call $A$-chords or $B$-chords
according to the splitting, as in figure 1.

\begin{center}
\Splittings\\
\medskip
{\bf Figure 1.}
\end{center}

 The Kauffman bracket polynomial
$\bd\
\in {\bf Z}[A^{\pm 1}]$ is defined by
\[\bd\ =\sum_{\mbox{states } s}\varphi_s,\] where
$\varphi_s=A^{a(s)-b(s)}(-A^2-A^{-2})^{|s|-1}$.

We write $\max(s)$ for the maximum degree of $\varphi_s$, so that
\[\max(s)=a(s)-b(s)+2|s|-2,\] and similarly
$\min(s)=a(s)-b(s)-2|s|+2$.

Suppose that a state $s'$ is given by changing $k$ of the
$A$-splittings of a state $s$ to $B$-splittings. It can be readily
shown that
\[\max(s')\le
\max(s),\] and that equality occurs if and only if $|s'|=|s|+k$.
Similarly, $\min(s')\le\min(s)$, with equality under the same
conditions.

  Write
$s_A$ and
$s_B$ for the {\em extreme states} in which all crossings are
labelled $A$ or all are labelled
$B$. It follows at once
\cite{K} that for any diagram $D$ with $c(D)$ crossings
\begin{eqnarray*}
\max\deg\bd&\le&\max(s_A)=c(D)+2|s_A|-2\\
\mbox{and }\min\deg\bd&\ge&\min(s_B)=-(c(D)+2|s_B|-2).
\end{eqnarray*}

We will call these bounds the {\em extreme degrees} for
$D$. They may of course be different from the actual minimum and
maximum degrees. Indeed both depend on the diagram chosen to
represent a given link, and they can be changed substantially by
Reidemeister moves, whereas $<D>$ is unaltered by moves II and III
and is simply multiplied by a power of $-A^3$ under move I.

The {\em spread} (or Laurent degree)
\[\beta(L)=\max\deg\bd-\min\deg\bd,\] which depends only on the Jones
polynomial of the link
$L$ represented by $D$, then satisfies the {\em extreme states bound}
\[\beta(L)\le 2c(D)+2(|s_A|+|s_B|-2).\]

We now give an algorithm for calculating the
coefficients of the terms of extreme degree, $\max(s_A)$ and
$\min(s_B)$, for a chosen diagram in terms of combinatorial features
of the state curves and splitting chords for each extreme state.
Thistlethwaite
\cite{Th} showed that for an {\em adequate} diagram both extreme
coefficients are
$\pm1$. Adequate diagrams include reduced
alternating diagrams, and so the bound shown gives an exact count of
the spread for any link with such a diagram. Thistlethwaite also found
a bound for the spread  expressed in terms of simple features of the
decomposition of a given diagram into maximal alternating tangles. We
show also how the non-alternating skeleton of the diagram can be
used to calculate the extreme states bound, which is in general stronger than Thistlethwaite's bound.

\section{The extreme coefficients}
In this section we give a combinatorial formula for the 
coefficient $a_{s_A}$ of the extreme term of degree $\max(s_A)$.

As noted in section 1, the only
contributions to $a_{s_A}$ come from states $s$ with $k$
$B$-splittings, such that
$|s|=|s_A|+k$, or equivalently $|s|=|s_A|+b(s)$. Each such state
$s$ then contributes
$(-1)^{|s_A|+k-1}$ to the extreme coefficient $a_{s_A}$.

Perform the $A$-splitting at each crossing of $D$ to get the
$|s_A|$ state circles for the extreme state $s_A$, marking
all the $A$-chords, as in figure 1.

Each state $s$ corresponds to a selection of a subset of $b(s) $ of
the $A$-chords at which the $B$-splitting is to be done in place of
the $A$-splitting. When a change of splitting takes place the number
of circles either increases or decreases by $1$, depending on whether
the ends of the splitting
$A$-chord lie on the same component or not. In order to finish with an
extra
$b(s)$ circles after performing $b(s)$ splittings we must increase the
number of circles at each splitting.

\begin{theorem}
Necessary and sufficient conditions for the splitting along a set of
$k$ $A$-chords of the state circles of $s_A$ to yield $k$ extra
curves are
\begin{itemize}
\item[(i)] the ends of each chord lie on the same 
state circle of $s_A$,
\item[(ii)] the ends of each pair of chords which lie on the same
circle do not  alternate in order round the circle.
\end{itemize}
\end{theorem}

\begin{Proof}
These conditions are clearly necessary, to ensure that the number of
circles increases after each splitting. They are also sufficient, by
induction on the number of chords, since after one splitting the
conditions are maintained. 
\end{Proof}

\begin{Def} A subset $C$ of the $A$-chords of the diagram $D$
is {\em independent} if it satisfies conditions (i) and (ii) above.
We include the case $C=\phi$.
\end{Def}

\begin{theorem} \label{maxdeg}
The coefficient of the term of degree $\max(s_A)$ is given by
\[a_{s_A}=(-1)^{|s_A|-1}\sum(-1)^{|C|},\] where the sum is taken over
all independent sets of $A$-chords $C$ in the state circles for $s_A$.
\end{theorem}

\begin{Proof} Each set of $A$-chords corresponds to a state of $D$, and
the set $C$ makes a contribution of $(-1)^{|s_A|+|C|-1}$ to the term
of degree $\max(s_A)$ if and only if $C$ is independent, as observed
above.
\end{Proof}

A diagram is {\em $+$ adequate} in  Thistlethwaite's sense when
$\max(s_A)>\max(s)$ for all states $s\neq s_A$. This is the case when
there are no non-empty independent sets of $A$-chords.

An exactly similar analysis of the term of degree $\min(s_B)$ can be
made in terms of the $B$-chords joining the state circles of the
state $s_B$. Here the coefficient is $(-1)^{|s_B|-1}\sum(-1)^{|C|}$,
with the sum taken over independent sets of $B$-chords. A diagram is
$-$ adequate when $\min(s_B)<\min(s)$ for all $s\neq s_B$, or
equivalently when there are no non-empty independent sets of
$B$-chords.

In these calculations the only chords that need to be considered are
those whose ends lie on the same state circle of
$s_A$ (or $s_B$). Indeed we can reformulate theorem \ref{maxdeg} in  a graph theory context as follows.

\begin{Def}
Let $K$ be a graph and let $C$ be a subset of the vertices of $K$.
Say that $C$ is {\em independent} if no two vertices of $C$ are
joined by an edge of $K$. Define an integer-valued  function $f$ on graphs
by
$f(K)=\sum (-1)^{|C|}$, where the sum is taken over all independent
subsets
$C$ of vertices of $K$, including the empty set. 
\end{Def}

Theorem \ref{maxdeg} can then be restated in terms of the function
$f$ for a suitable graph.
\begin{theorem}[graphical version]
Let $K$ be  the
`intersection graph', in the sense of Lando \cite{Lando}, of
the $A$-chords with endpoints on the same state circle for $s_A$,
namely   the graph
whose vertices are these chords, 
with an edge joining each pair of chords whose ends occur
alternately on the same circle. Then  the required coefficient
$a_{s_A}$ is just $(-1)^{|s_A|-1}f(K)$. 
\end{theorem}

Calculation of $f$  for a general graph $K$ is simplified by the use of some readily
established properties.
\begin{property}
 [Recursion]

Let $K-v$ be the subgraph of $K$ given by
deleting a vertex $v$ and its incident edges, and let $K-N(v)$ be
given by deleting the immediate neighbours of $v$ along with their
incident edges. Then
\[f(K)=f(K-v)-f(K-N(v)).\] 
\end{property}
\begin{Proof}
This follows by grouping the subsets $C$ into those which do and
those which do not contain the vertex $v$.
\end{Proof}

\begin{property}[Multiplication under disjoint union] \label{multiplication}

Let
$K=K_1\cup K_2$ with
$K_1\cap K_2=\phi$ then $f(K)=f(K_1)f(K_2)$.
\end{property}

\begin{Proof}
This follows by induction on the number of vertices in $K_2$, and
recursion.
\end{Proof}

\begin{property} [Duplication]

 If an extra vertex $v$ is inserted which is {\em
not} joined to one vertex $w$ but which is joined to {\em all} the
neighbours of
$w$ (and possibly to other vertices as well)  then $f$ is unchanged.
\end{property}

\begin{Proof}
This follows from the first two properties, noting that here $K-N(v)$
has an isolated vertex $w$ and that $f=0$  on a graph with a
single vertex.
\end{Proof}

\begin{remark} The intersection graphs which arise here are naturally
bipartite graphs, as the planar diagram of chords has one set of
non-intersecting chords inside each circle, and another set of
non-intersecting chords outside the curvecircle. Intersections in
the graph can only take place between the two different types of chord,
and can be realised by redrawing so that both sets of chords lie
inside the circle, when the condition that the endpoints occur
alternately will correspond to an intersection of a pair of chords.
Such graphs are sometimes known as {\em circle graphs}.
\end{remark}

 While the calculation of $f(K)$
clearly depends only on the intersection matrix of the graph $K$, we
do not have a simple formula for
$f$ in terms of this matrix. Clearly if there is just a single chord
then $f=0$, and equally if there is no intersection among any of the
chords then again $f=0$ by the multiplicative property. On the other
hand there are plenty of examples where the graph is non-empty and the
value of $f$ is non-zero, giving us exact bounds on the spread of the
bracket polynomial beyond the cases of
$\pm$ adequate diagrams.

\section{Alternating tangle decompositions}
In using theorem \ref{maxdeg} to calculate $a_{s_A}$ it is  enough
to consider the state circles individually, because of the
multiplicative property \ref{multiplication} of $f$.

Where a diagram has a substantial number of alternating edges
there will in general be many of the extreme state circles with no
cross-chords. These can therefore be ignored completely in the
calculation.

 Consider the projection of the diagram
$D$ as a
$4$-valent planar graph, which we call the {\em projection graph} of $D$.  Each
edge is either {\em alternating} or {\em non-alternating} according
to the crossings at its ends. The non-alternating edges are of two
types, {\em over} and {\em under}, indicated by $+$ and $-$
respectively. 

  The state circles of
$s_A$ and the
$A$-chords are constructed by separating the vertices of the projection graph
 slightly and inserting the appropriate chord. They can be generated
dynamically as circuits in the projection graph by turning right at each undercrossing of $D$, and left at each overcrossing. We assume throughout that
$D$ is not a split diagram, and is {\em reduced} in the sense that it has
no cut-vertex. The closure of each complementary region of its graph
in
$S^2$ is then a disc.  Any state circle for
$s_A$ consisting entirely of alternating edges forms the boundary of
one of these discs, as the remainder of the graph lies entirely on the same side of the state circle which consequently has no cross-chords. In using
theorem \ref{maxdeg} we need then only consider chords with ends on
those state circles which include some non-alternating edges.

To find these systematically we draw the graph $G$ which is dual to the
non-alternating edges in the projection graph of $D$; there is one vertex of $G$ in each complementary disc of the projection graph whose boundary contains non-alternating edges.

When we superimpose $G$ on the original knot diagram $D$ we find that
the complementary regions of $G$ are discs  if $G$ is connected, or
more generally discs with holes, which meet $D$ in alternating
tangles. The intersections of $D$ with these complementary regions
are the maximal alternating pieces as defined by Thistlethwaite \cite{Th};
 when $G$ is not
connected some of them will be tangles in a disc with holes, rather
than a classical tangle in a disc. In our setting, Thistlethwaite's `channels' are the planar neighbourhoods of the components of $G$.

In any event, the major part of each tangle can be ignored in
making our calculations, and we  concentrate on the graph $G$,
which we call the {\em non-alternating skeleton} of $D$.

Each state circle for $s_A$ made up of alternating edges bounds a
disc lying entirely in one of the alternating tangles. Only the state circles with
non-alternating edges intersect the skeleton $G$; we show how to
recover them up to isotopy by making a standard splitting of the
graph $G$. We then add the `boundary information' about the $A$-chords with ends
on these circles to complete the data needed in the calculations of
theorem \ref{maxdeg}.

\section{The non-alternating extreme state circles}

The non-alternating edges of the diagram $D$ and hence the edges of its non-alternating skeleton $G$ come in two types, {\em over}, labelled $+$ and {\em under}, labelled $-$.  As we trace out a state circle of $s_A$ which contains some non-alternating edges we will come to a non-alternating over edge, labelled $+$, where the circle will cross $G$. It will then continue past some crossings, turning left  each time, until it reaches the next non-alternating edge, necessarily an under edge, with sign $-$, where it again crosses $G$. These two edges of $G$  have a common vertex, lying in the same complementary region of the projection graph, and so the local picture of the projection graph and $G$ will look like figure 2.

\begin{center}
\figtwo\\
{\bf Figure 2}
\end{center}
\medskip

The two edges of $G$ are then isotopic, in the complement of $G\cup D$ to this segment of the $s_A$ state circle.

The $+$ and $-$ edges must alternate around each vertex of $G$, and all pairs will correspond in this way to pieces of states circles. So when we break the graph $G$ apart at each vertex by pairing adjacent $+$ and $-$ edges, matching each $+$ edge with the next $-$ edge in the anticlockwise sense, as in figure 3,

\begin{center}
$G\ =\ $\figthreea\ \quad $\to$ \quad $G_A\ =\ $\ \figthreeb\\
{\bf Figure 3}
\end{center}
\medskip
the resulting curves, which we denote by $G_A$, are isotopic to the $s_A$ curves with non-alternating edges. These are the only $s_A$ circles which can appear in our extreme term calculation. 

Separating all the vertices of $G$ in the opposite sense will similarly yield curves $G_B$ isotopic to  the non-alternating $s_B$ circles.

\subsection{Boundary information}

Having found all the $s_A$ circles  needed for theorem \ref{maxdeg} we can identify those $A$-chords which may be involved in the formula. Recalling that only $A$-chords with both ends on the same circle will contribute, we may treat the  components of $G$ separately, and combine the results by use of the multiplicative property of $f$.  

From our picture of the construction of the non-alternating $s_A$ circles we see that the possible $A$-chords occur in a complementary region of $G$ where there is an arc across the region which passes through just one crossing of the projection graph, as in figure 4. 
\begin{center}
\figfour\\
{\bf Figure 4}
\end{center}

 To give an $A$-chord the crossing must be approached through the $B$-quadrants, in terms of the original diagram $D$. The $A$-chord which arises from the crossing is then isotopic to this arc when the $s_A$ circle is isotoped to $G_A$.
  The union of such arcs drawn across the complementary regions of $G$ is the `boundary information' which we need for the $A$-chords, with a similar union of arcs drawn for the $B$-chords.

\subsection{Algorithm for finding the extreme coefficients}

Given a diagram $D$ we can assemble the data needed to apply theorem \ref{maxdeg} as follows.

\noindent
{\bf Step 1}. Construct the non-alternating skeleton $G$.

\noindent
{\bf Step 2}. Add the boundary information.

For each component of $G$ consider separately each tangle defined by $D$ in the complementary regions of this component.  Draw any arcs across each tangle which meet $D$ in just one  crossing approached through the $B$-quadrants. In general there will be a relatively small number of these, as they can only involve the complementary regions of $D$ which are adjacent to the boundary of the tangle, hence the term `boundary information'. The result is to decorate $G$ with a number of non-intersecting chords drawn across the complementary regions.

\noindent
{\bf
Step 3}.
Separate $G$ into the non-alternating $s_A$ circles $G_A$ by splitting apart at the vertices, as above, while retaining the decorating chords.

\noindent
{\bf
Step 4}. For each circle separately calculate the value of $f$ on the intersection graph given by the cross-chords. Multiply the values, to give the eventual value of the coefficient $a_{s_A}$ of top degree.

Repeat steps 2-4 with the $B$-chords across the tangles, and the splitting of $G$ into $s_B$ circles, to find the lowest degree coefficient similarly.
\begin{Proof}
Separating $G$ gives the circles $G_A$ which are isotopic to the non-alternating $s_A$ state circles. The boundary information recovers all possible $A$-chords with ends on the same circle. The calculations then follow from theorem \ref{maxdeg} and property \ref{multiplication} of $f$.
\end{Proof}

\section{States surfaces}

One of the neatest techniques in the proof of the original results
about alternating diagrams is the use of  `states surfaces'.  Each
state $s$ of a diagram $D$ has a dual state $\hat s$ defined by
changing all the markers of $s$. In particular the extreme states
$s_A$ and $s_B$ are dual to each other.  The states surface for $s$
is a closed orientable surface with Euler characteristic
$|s|+|\hat s|-c(D)$. The {\em extreme states surface}
$F$ for the states $s_A$ or $s_B$ then has Euler characteristic
$\chi(F)=|s_A|+|s_B|-c(D)$ and the extreme states bound for the spread
of $\bd$, which is
$2(c(D)+|s_A|+|s_B|-2)$, can  be written as
$4c(D)-4g(F)$ in terms of the genus $g(F)$ of 
$F$.

In \cite{T} Turaev gives a construction for the extreme states surface
$F$ in which  discs round each crossing of $D$ are connected by an
untwisted band for each alternating edge, and a half-twisted band for
each non-alternating edge. This gives a surface with
$|s_A|+|s_B|$ boundary components, and yields $F$ when they are
capped off by discs. 

Make this construction with the non-alternating  skeleton $G$ in place, inserting first only the bands for the alternating edges. The boundary of the resulting planar surface  includes all state circles for $s_A$ and $s_B$ made of alternating edges only. These all lie in complementary regions of $G$, along with circles parallel to the boundary of each complementary region.  Capping off the alternating state circles then gives the complement of a neighbourhood of $G$.

The surface $F$ is completed by adding a twisted band across each edge of $G$ and capping off the boundary of the resulting surface.  The boundary curves of this surface can be readily identified with the non-alternating state circles given by separating the vertices of $G$ to yield the curves $G_A$ and $G_B$.
Then
\begin{eqnarray*} \chi(F)&=&2-\chi(G)-e(G)+|G_A|+|G_B|\\
&=&2-v(G)+|G_A|+|G_B|,
\end{eqnarray*} where $G$
has $v(G)$ vertices and
$e(G)$ edges. The extreme states bound is then
\begin{eqnarray*} 4c(D)-4g(F)&=&4c(D)+2\chi(F)-4\\
&=&4c(D)+2(|G_A|+|G_B|)-2v(G).
\end{eqnarray*}

The extreme states bound for the spread can thus be found readily in
terms of the non-alternating skeleton $G$, as an embedded graph (so as to find $G_A$ and $G_B$).

\begin{theorem} The extreme states bound is lower than Thistlethwaite's bound in general.
\end{theorem}
\begin{Proof} Thistlethwaite's bound for the spread $\max\deg\bd-\min\deg\bd$ is given in terms of the number of alternating tangles $n$, (the number of complementary regions of $G$), and the number of non-alternating edges $\nu$ (=$e(G)$).
Explicitly, his bound is
$ 4c(D)+4(n-1)-2\nu$.

Suppose that $G$ has $r$ components, so that  $n=r+1-v(G)+e(G)$.
Now construct a surface from $G$ by putting a disc at each vertex, and joining them by a twisted band for each edge. The boundary can again be regarded as the curves $G_A$ and $G_B$. Cap these off to give a closed surface $F^*$ with $r$ components, so that $\chi(F^*)\le 2r$. Then\[0\le4r-2\chi(F^*)=4r-2(v(G)-e(G)+|G_A|+|G_B|),\] and so the extreme states bound of $4c(D)-4g(F)$ satisfies
\begin{eqnarray*}4c(D)-4g(F)&=&4c(D)+2(|G_A|+|G_B|)-2v(G)\\
&\le&4c(D)+4r-4v(G)+4e(G)-2e(G)\\
&=& 4c(D)+4(n-1)-2\nu.
\end{eqnarray*}
\end{Proof}

\section{Some examples}

In figure 5 we show a diagram with its non-alternating skeleton $G$, and the two sets of curves $G_A$ and $G_B$ resulting from  splitting $G$. In this case $e(G)=10$ and $\chi(G)=-2$, while $|G_A|=|G_B|=2$. The extreme states bound  for the spread of the bracket polynomial is then $4c+2\times4-2(\chi(G)+e(G))=4c-8$, giving a bound of $c-2=19$ for the spread of its Jones polynomial.

\begin{center}
\figfive\\
\medskip
$G_A\ =\ $ \figfivega\ ,\qquad $G_B\ =\ $ \figfivegb\ .\\
{\bf Figure 5}
\end{center}

In figure 6 we show separately the reducing $A$-chords and $B$-chords on $G$. After splitting $G$ and retaining only the cross-chords with both ends on the same circle we get the essential boundary information shown in figure 7.

\begin{center}
$A$-chords,\ \figsixa\ ,\qquad $B$-chords,\ \figsixb\ .\\
{\bf Figure 6}
\end{center}
\begin{center}
$G_A =\ $\figsevena\ ,$\quad G_B\ =\ $\figsevenb\ .\\
\medskip
 {\bf Figure 7.}
\end{center}

Since there are no chords on $G_B$ the diagram is $-$ adequate, but in view of the three chords on $G_A$ it is not $+$ adequate. Apply the function $f$ to the circle graph of $G_A$ to calculate $\hat a_{s_A}=(-1)^{|s_A|-1}a_{s_A}=-1$. Both extreme coefficients are then non-zero and we deduce
that the exact spread of the Jones polynomial is $c-2$.

The diagram in figure 8 has the same non-alternating skeleton $G$ and only differs from figure 5 in that the alternating tangle in one of the complementary regions of $G$ has been rotated.

\begin{center}
\figeight\\
{\bf Figure 8}
\end{center}

Thus $G_A$ and $G_B$, and the bound of $c-2$ for the spread of the Jones polynomial are unaltered. However the two reducing chords in the rotated tangle now lie in a different way relative to $G$, as shown in figure 9, and the new boundary information for $G_A$ is shown in figure 10.

\begin{center}
\fignine\\
{\bf Figure 9}
\end{center}

\begin{center}
$G_A\ =\ $\figten\ .\\
\medskip
 {\bf Figure 10.}
\end{center} 

The diagram is still $-$ adequate, but $f=0$ for the
circle graph of $G_A$ and so
$a_{s_A}=0$. The spread of the Jones polynomial is then at most $c-3=18$.

\subsection{Calculations on Rolfsen's tables}
Table 1 shows the reducing chords, and the
corresponding coefficients, for the diagrams of knots up
to 10 crossings. The source of the table  is
Rolfsen's knot diagram table of 10 crossings
or less. Since the maximal and minimal terms
of the alternating knots are already known, we
list data for non-alternating knots only.
The non-alternating skeleton for all these diagrams in Rolfsen's table consists of a single circle, so $G_A$ and $G_B$ are a single curve in each case.

In each row of the table we show the diagrams of the
$A$-reducing and $B$-reducing chords. The extreme coefficients are
given by  $a_{s_A}=(-1)^{\abs{s_A}-1}\hat a_{s_A}$ and
$a_{s_B}=(-1)^{\abs{s_B}-1}\hat a_{s_B}$, where $\hat a_{s_A}$ and $\hat
a_{s_B}$ are calculated directly from the intersection graph using the
function $f$.  The number 
$\hat\beta$ is the extreme
states  bound for the spread of the Jones polynomial.  This is equal to
the actual spread$\beta$ when
$a_{s_A}\not= 0$ and $a_{s_B}\not= 0$. The value of $\beta$ in
other cases is noted in the final column of the table for comparison,
calculated directly from the Jones polynomial.

Where the diagram of $G_A$ or $G_B$ has no chords the knot diagram
is  
$\pm$ adequate, so we can see that all knots of 10 crossings or less
are
$+$ adequate or $-$ adequate. In particular, $10_{152}, 10_{153},
10_{154}$ are the only adequate knots of 10 crossings or less which
are non-alternating.

The two different values in $\hat a_{s_B}$ for $10_{144}$ reflect
a mistake in the diagram in Rolfsen's table. For the knot which is
presented by Rolfsen's diagram, $\hat a_{s_B}=-1$ while for the genuine
knot
$10_{144}$, $\hat a_{s_B}=2.$
\bigskip

\begin{center}
Table 1. Reducing chords and extreme coefficients up to 10 crossings.
\\[5pt]
\begin{tabular}{|c||c|c|c|c|c|c|c|c|} \hline
  Knot& \ \ $\hat\beta$\ \    &$\hat a_{s_B}$  &   chords on $G_B$   & $\hat a_{s_A}$ &chords on $G_A$ &
$\beta(\seq{L})$\\ \hline\hline
  $8_{19}$    &    6    &       0  & \input{8_19.lp}     &       1  & 
\input{trivial.lp}  &5 \\ \hline
  $8_{20}$    &    6    &       1  & \input{8_20.lp}     &       1  & 
\input{trivial.lp}   &\\ \hline
  $8_{21}$    &    6    &       2  & \input{8_21.lp}     &       1  & 
\input{trivial.lp}   &\\ \hline
  $9_{42}$    &    8    &       0  & \input{9_42.lp}     &       1  & 
\input{trivial.lp}   &6\\ \hline
  $9_{43}$    &    8    &       0  & \input{9_43.lp}     &       1  & 
\input{trivial.lp}   &7\\ \hline
  $9_{44}$    &    7    &       1  & \input{9_44.lp}     &       1  & 
\input{trivial.lp}   &\\ \hline
  $9_{45}$    &    8    &       0  & \input{9_45.lp}     &       1  & 
\input{trivial.lp}   &7\\ \hline
  $9_{46}$    &    8    &       0  & \input{9_45.lp}     &       1  & 
\input{trivial.lp}   &6\\ \hline
  $9_{47}$    &    8    &       0  & \input{9_42.lp}     &       1  & 
\input{trivial.lp}   &7\\ \hline
  $9_{48}$    &    7    &       2  & \input{9_48.lp}     &       1  
&\input{trivial.lp}   & \\ \hline
  $9_{49}$    &    8    &       0  & \input{9_42.lp}     &       1  & 
\input{trivial.lp}   &7\\ \hline
 $10_{124}$    &    9    &       0  & \input{9_42.lp}   &       1  & 
\input{trivial.lp}   &6\\ \hline
  $10_{125}$    &    9    &       0  & \input{9_45.lp}   &       1  & 
\input{trivial.lp}   &8\\ \hline
  $10_{126}$    &    9    &       0  & \input{10_126.lp}   &       1 & 
\input{trivial.lp}   &8  \\ \hline
  $10_{127}$    &    9    &       0  & \input{9_42.lp}     &       1  & 
\input{trivial.lp}   &8 \\ \hline
  $10_{128}$    &    8    &       0  & \input{10_128.lp}   &       1  & 
\input{trivial.lp}   &6 \\ \hline
  $10_{129}$    &    9    &       0  & \input{9_42.lp}     &       1  & 
\input{trivial.lp}   &8 \\ \hline
 $10_{130}$    &    9    &       0  & \input{9_42.lp}     &       1  & 
\input{trivial.lp}    &8 \\ \hline
  $10_{131}$    &    9    &       0  & \input{9_45.lp}     &       1  & 
\input{trivial.lp}    &8 \\ \hline
  $10_{132}$    &    7    &       0  & \input{10_132.lp}     &       1  &  
\input{trivial.lp}    &5 \\ \hline
    $10_{133}$    &    8    &       1  & \input{10_133.lp}     &       1  & 
\input{trivial.lp}    & \\ \hline
    $10_{134}$    &    8    &       1  & \input{10_141.lp}     &       1  & 
\input{trivial.lp}     &\\ \hline
  $10_{135}$    &    9    &       0  & \input{9_42.lp}     &       1  &  
\input{trivial.lp}  &8\\ \hline
  $10_{136}$    &    9    &       0  & \input{9_43.lp}     &       1 & 
\input{trivial.lp}    &7 \\ \hline
  $10_{137}$    &    9    &       0  & \input{9_43.lp}     &       1 & 
\input{trivial.lp}    &8 \\ \hline
  $10_{138}$    &    9    &       0  & \input{9_42.lp}     &       1  & 
\input{trivial.lp}    &8 \\ \hline
  $10_{139}$    &    8    &       -1  & \input{10_139.lp}    &       1  & 
\input{trivial.lp}    & \\ \hline
   $10_{140}$    &    9    &       0  & \input{9_45.lp}     &       1  & 
\input{trivial.lp}    &7 \\ \hline
   $10_{141}$    &    8    &       1  & \input{10_141.lp}     &       1  & 
\input{trivial.lp}    & \\ \hline
  $10_{142}$    &    9    &       0  & \input{10_142.lp}     &       1  & 
\input{trivial.lp}    &6 \\ \hline
  $10_{143}$    &    9    &       0  & \input{9_42.lp}     &       1  &  
\input{trivial.lp}    &8 \\ \hline
  $10_{144}$    &    8    &       -1(2)&\input{10_144.lp}   &       1  & 
\input{trivial.lp}    & \\ \hline
\end{tabular}\\
\end{center}
\vfill\eject
\begin{center}
Table 1. Continued.\\[5pt] 
\begin{tabular}{|c||c|c|c|c|c|c|c|c|} \hline
  Knot& \ \ $\hat\beta$\ \    &$\hat a_{s_B}$  &  chords  on $G_B$                 & $\hat 
a_{s_A}$ & chords on $G_A$ & $\beta(\seq{L})$\\ \hline\hline
  $10_{145}$    &    8    &       -1  & \input{10_139.lp}    &       1  & 
\input{trivial.lp}    & \\ \hline
$10_{146}$    & 9       &       0  & \input{9_45.lp}     & 1         & 
\input{trivial.lp}   &8\\ \hline
  $10_{147}$    &    8    &       1  & \input{10_141.lp}     &       1  & 
\input{trivial.lp}    & \\ \hline
  $10_{148}$    &    9    &       0  & \input{9_42.lp}     &       1  & 
\input{trivial.lp}     &8\\ \hline
  $10_{149}$    &    9    &       0  & \input{9_45.lp}     &       1  & 
\input{trivial.lp}     &8\\ \hline
 $10_{150}$    &    8    &       1  & \input{10_141.lp}     &       1  & 
\input{trivial.lp}    & \\ \hline
  $10_{151}$    &    9    &       0  & \input{9_42.lp}     &    1  &  
\input{trivial.lp}    &8 \\ \hline
  $10_{152}$    &    9    &       1  & \input{trivial.lp}     &       1  & 
\input{trivial.lp}    & \\ \hline
  $10_{153}$    &    9    &       1  & \input{trivial.lp}     &       1  & 
\input{trivial.lp}    & \\ \hline
  $10_{154}$    &    9    &       1  & \input{trivial.lp}     &       1  & 
\input{trivial.lp}    & \\ \hline
  $10_{155}$    &    9    &       1  & \input{trivial.lp}     &   0  &  
\input{9_45.lp}  &8\\ \hline
  $10_{156}$    &    8    &       1  & \input{trivial.lp}      &       1  & 
\input{10_156.lp}  & \\ \hline
  $10_{157}$    &    8    &       1  & \input{trivial.lp}      &       2  & 
\input{10_157.lp}   &\\ \hline
  $10_{158}$    &    9    &       0  & \input{9_45.lp}     &       1  & 
\input{trivial.lp}   &8  \\ \hline
  $10_{159}$    &    9    &       0  & \input{9_45.lp}     &       1  & 
\input{trivial.lp}     &8\\ \hline
  $10_{160}$    &    8    &       0  & \input{10_160.lp}     &       1  & 
\input{trivial.lp}    &7 \\ \hline
  $10_{161}$    &    8    &       -1 & \input{10_161.lp}     &       1  & 
\input{trivial.lp}    & \\ \hline
  $10_{162}$    &    8    &       -1 & \input{10_161.lp}       &       1  & 
\input{trivial.lp}     &\\ \hline
  $10_{163}$    &    8    &       2  & \input{10_163.lp}     &       1  & 
\input{trivial.lp}    & \\ \hline
  $10_{164}$    &    9    &       0  & \input{9_45.lp}     &       1  & 
\input{trivial.lp}    &8 \\ \hline
  $10_{165}$    &    8    &       2  & \input{10_165.lp}     &       1 & 
\input{trivial.lp}     &\\ \hline
  $10_{166}$    &    9    &       0  & \input{9_42.lp}     &       1  & 
\input{trivial.lp}    &8 \\ \hline
\end{tabular}
\end{center}
\medskip

\subsection{Realisation of extreme coefficients}
We finish with some results about the range of possible values of the extreme coefficients.  It is
certainly possible to find a graph $K$ with $f(K)=n$ for any chosen integer $n$. Indeed
it is easy to see that $f(K_{n+1})=-n$, where $K_{n+1}$ is the
complete graph on
$n+1$ vertices. On the other hand the circle graphs which determine
the extreme coefficients form a proper subset of all bipartite
graphs, and $f(K_{m,n})=-1$ for the complete bipartite graph
$K_{m,n}$.

We initially wondered whether any values of $f$ besides $0$ and $\pm2^k$ were possible for the
extreme coefficients. We then managed to find a circle graph with $f=3$,
 illustrated in figure 11 along with its realisation by chords,
and used it to produce a link with extreme coefficient $3$.

 \begin{center}
\Doublehex\quad\qquad\Doublehexcircle\\
\medskip
{\bf Figure 11.}
\end{center}

This can be done by
replacing each chord, $\Chord$, in the set $E\cup F$ of chords in the
realisation, by a single crossing $\Ychord$. The circle graph of
figure 11 produces in this way  the link
 shown in figure 12, whose
bracket polynomial is $3A^{13}-2A^9+4A^5-A+4A^{-3}-A^{-7}+A^{-11}$.

 \begin{center}
 \Doublehexknot\\
\medskip
{\bf Figure 12.}
\end{center}

More generally, given any set $E$ of non-intersecting chords inside a circle,
and another set $F$ of non-intersecting chords outside the same circle the same
procedure will construct a knot having a single curve $G_A$ with $E$ and
$F$ as its
$A$-reducing chords,  hence with one extreme coefficient given by the
intersection graph of these chords. As in the example above, the other extreme
coefficient may be zero. A more elaborate construction, using for
example $\Ytwochord$ in place of some or all of the chords, can be made to
ensure that the tangles used in the construction are
$B$-reduced, while retaining the same or parallel families of
$A$-reducing chords. Starting from a circle graph with value $f$ this will lead
to a knot or link with a 
$-$adequate diagram, whose extreme coefficients are then $f$ and $1$, up to
sign. It is equally easy to extend this so that both extreme coefficients are
any chosen values of $f$ for a circle graph. 

 The coefficients of the maximal and minimal degree terms in the
Jones polynomial may not in general be values of $f$, since the
knot may not have a diagram for which they appear as the extreme
coefficients. Since the original version of this paper was written
there has been a nice construction due to Manchon \cite{Manchon}
giving  circle
graphs (or equivalently the families of chords
$E$ and
$F$) which realise every integer  value of
$f$.

\medskip

\begin{leftline}
\noindent{\ }\\
Department of Mathematics,\\
College of Natural Sciences,\\
Kyungpook National University,\\
Taegu, 702-701, Korea.\\

\noindent
Department of Mathematical Sciences,\\ 
University of Liverpool,\\
Peach St.\\
Liverpool, L69 7ZL, England.
\end{leftline}

\medskip

\noindent
email {\tt ybae@knu.ac.kr  morton@liv.ac.uk}

\end{document}